\documentclass[a4paper,12pt]{article}
\usepackage{amssymb}
\usepackage{amsfonts}
\usepackage{amsmath}
\usepackage[top=2cm,bottom=2cm,right=2cm,left=2cm]{geometry}
\usepackage[backend=bibtex]{biblatex}

\newcommand{\be}{\begin{equation}}
\newcommand{\ee}{\end{equation}}
\newcommand{\op}{ \oplus_q  }
\newcommand{\om}{ \ominus_q  }

\newcommand{\bea}{\begin{eqnarray}}
\newcommand{\eea}{\end{eqnarray}}
\begin{document}

\title{Connection formulae for new $q$-deformed Laguerre-Gould-Hopper
polynomials }
\author{Sama Arjika$^{\natural ,\ddag }$ and Zouha\"{\i}r Mouayn$^{\ast }$}
\date{$^{\natural }${\footnotesize \ Department of Mathematics and Computer
Sciences, Faculty of Sciences, \vspace{-0.5em}\\ University of Agadez, BP. 199, Agadez, Niger}\\
$^{\ddag }${\footnotesize \ International Chair in Mathematical Physics and
Applications {\scriptsize(ICMPA-UNESCO Chair)}\vspace{-0.5em}\\
 University of Abomey-Calavi, 072,  BP. 50,
Cotonou,  Benin}\medskip\\
$^{\ast }${\footnotesize \ Department of Mathematics, Faculty of Sciences
and Technics (M'Ghila)\\ BP. 523, B\'{e}ni Mellal, Morocco}}
\maketitle

\begin{abstract}
We introduce new families of  $q$%
-deformed 2D Laguerre-Gould-Hopper polynomials. For  these  polynomials  we establish connection formulae which extend some known
ones.
\end{abstract}

\section{Introduction}
The $q$-analysis  can be traced back to the earlier works of L. J. Rogers  [1]. It has wideranging applications in the analytic number theory  and $q$-deformation of well-known functions [2] as well as in the  study of solvable models in statistical mechanics [3]. During the 80's the interest on this  analysis increased with quantum groups theory  with which models of $q$-deformed oscillators have been developed [4]. The $q$-analogs of boson operators have been defined in [5] where the corresponding wavefunctions were constructed  in terms of the continuous $q$-Hermite polynomials of Rogers and other polynomials. Actually, known models of $q$-oscillators are closely related with $q$-orthogonal polynomials.

 The  connection formulae for  orthogonal polynomials are useful in mathematical analysis and also have applications in quantum mechanics, such as finding the relationships between the wavefunctions of some  potential functions used for describing physical and chemical properties in atomic and molecular systems [6].
 
 Here we deal with connection formulae for a new class of $q$-deformed Laguerre-Gould-Hopper polynomials which we introduce in an operatorial  way  by acting on certain newly defined $q$-deformed Gould-Hopper generalized Hermite polynomials [7]. We also consider particular cases of these formulae, which enables  us to recover known ones in the literature.
 
 The paper is organized as follows. In section 2, we prepare some needed notations and definitions. We also introduce a family of $q$-deformed Laguerre polynomials and express them by a  generating function. In section 3, we define a new class of $q$-deformed Laguerre-Gould-Hopper polynomials for which we give a generating function and we establish some connection formulae. Section 4 is devoted to discuss some particular cases of the obtained results. 

\newpage
\section{Notations and definitions}

Here, we list notations and special functions we will be using. For the relevant properties and
definitions we refer to [8-10]. We also introduce a family of $q$-deformed 2D Laguerre polynomials and we write down their generating function.
\begin{enumerate}

\item The $q$-analogues ($ \:q\in\mathbb{C}\setminus\{1\}$) of a natural  number $n$,  the factorial and  semifactorial functions are defined by  
\begin{equation}
\lbrack n]_{q}\equiv \sum_{k=1}^{n}q^{k-1},\: [n]_q!=\displaystyle\prod_{k=1}^{n}[k]_q!,\: \lbrack 0]_{q}!\equiv 1,  \tag{2.1}
\end{equation}%
\begin{equation*}
 \left[ mk\right] _q!!:=\prod\limits_{l=1}^{k}\left[ ml\right] _{q},\qquad \left[ 0\right] !!:=1. \tag{2.2}
\end{equation*}
\item The Euler-Heine-Jackson $q$-difference operator is defined by 
\begin{equation}
D_{x}^{q}[f](x):=\frac{f(x)-f(qx)}{(1-q)x},\qquad x\neq 0.  \tag{2.3}
\end{equation}%
Its inverse $(D_{x}^{q})^{-1}$ is defined in such a way that

\begin{equation}
(D_{x}^{q})^{-n}(1)=\frac{x^{n}}{[n]_{q}!}.  \tag{2.4}
\end{equation}%
\item The Gauss $q$-binomial coefficient is  given by%
\begin{equation}
\begin{bmatrix} n\\ k \end{bmatrix}_q:=\frac{[n]_{q}!}{[n-k]_{q}![k]_{q}!}=\frac{(q;q)_{n}}{%
(q;q)_{n-k}(q;q)_{k}},\text{ } \: k=0,1,\cdots,n.  \tag{2.5}
\end{equation}%
where 
\begin{equation}
(a;q)_{n}:=\medskip \prod\limits_{k=0}^{n-1}\left( 1-aq^{k}\right)
,(a;q)_{0}=1  \tag{2.6}
\end{equation}%
\ denotes the $q$-shifted factorial.

\item The Jackson-Hahn-Cigler (JHC)  $q$-addition is the function 
\begin{equation*}
(x\oplus _{q}y)^{n}:=\sum_{k=0}^{n}\begin{bmatrix} n\\ k \end{bmatrix}_q q^{\binom{k}{2}
}x^{n-k}y^{k}
\end{equation*}
\begin{equation*}
\qquad\qquad\qquad\qquad\qquad\qquad\qquad\qquad\:\:\;\:=x^n\left(\frac{-y}{x};q\right)_n\equiv P_{n,q}(x,y),\: n=0,1,2,\cdots,\tag{2.7}
\end{equation*}
The JHC $q$-substruction is defined by 
\begin{equation*}
(x\ominus _{q}y)^n:=P_{n,q}(x,-y),\: n=0,1,\cdots,  
\end{equation*}

\item Let $\ D_{\rho }=\{\xi \in \mathbb{C},|\xi |<\rho \}$ the
complex disk of radius $\rho >0$ and  $ F(z):=\sum_{n=0}^{\infty }c_{n}z^{n}$, $z\ in D_{\rho }$. Define
the formal series 
\begin{equation}
F(x\oplus _{q}y)=\sum_{n=0}^{\infty }c_{n}(x\oplus _{q}y)^{n}.  \tag{2.8}
\end{equation}

\item For $m=0,1,2,...$, let $e_{q}$ and $E_{q^{m}}$ denote the
exponential functions defined by 
\begin{equation}
e_{q}(x):=\sum_{n=0}^{\infty }\frac{x^{n}}{[n]_{q}!}  \tag{2.9}
\end{equation}%
and 
\begin{equation}
E_{q^{m}}\left( x\right) :=\sum_{n=0}^{\infty }\frac{q^{m\binom{n}{2}}}{%
[n]_{q^{m}}!}x^{n}.  \tag{2.10}
\end{equation}%
These functions satisfy $e_{q}(x)E_{q^{m}}(y)=e_{q}(x\ominus _{q,q^{m}}y)$
where 
\begin{equation}
(a\ominus _{q,q^{m}}b)^{n}=[n]_{q}!\sum_{k=0}^{n}\frac{(-1)^{k}q^{m\binom{k}{%
2}}}{[n-k]_{q}![k]_{q^{m}}!}a^{n-k}b^{k}  \tag{2.11}
\end{equation}%
with $(a\ominus _{q,q^{m}}b)^{0}:=1.$ For $m=1$, we have the rules 
\begin{equation*}
e_{q}(a)E_{q}(b)=e_{q}(a\oplus _{q}b),\quad e_{q}(a)E_{q}(-a)=1. \tag{2.12}
\end{equation*}

\item For $m=0,1,2,\cdots ,$ the function  
\begin{equation}
\varepsilon _{n,q}^{m}(x):=\sum_{n=0}^{\infty }\frac{(-1)^{k}q^{m\binom{k}{2}%
}}{[k]_{q^{m}}!\left[ n+k\right] _{q^{m}}!}x^{k},\quad n=0,1,2,\cdots  \tag{2.13}
\end{equation}%
denotes the $n$th order of a $(q,m)$-deformed Bessel-Tricomi function with $\varepsilon _{0,q}^{m}(x):=\varepsilon _{q}^{m}(x)$.

\item  The  Gould-Hopper generalized polynomials are defined as ([11], p.58): 
\begin{equation}
g_n^m(x,y):=n!\displaystyle\sum_{k=0}^{\lfloor \frac{n}{m}\rfloor} \frac{x^{n-mk}y^k}{k!(n-mk)!} \tag{2.14} 
\end{equation}
where $\left\lfloor a\right\rfloor $ stands
for the greatest integer not exceeding $a.$

\item We define  a class of $q$-deformed two variables (2D)
Laguerre polynomials  ($q$-2DLP) as 
\begin{equation}
_{m}L_{n}(x,y|q):=E_{q^{m}}\left( (D_{x}^{q^{m}})^{-1}(D_{y}^{q})^{m}\right)
(y^{n}),\:m,n=0,1,2,\cdots.  \tag{2.15}
\end{equation}%
Explicitly,
\begin{equation}
_{m}L_{n}(x,y|q):=[n]_{q}!\sum_{k=0}^{\left\lfloor \frac{n}{m}\right\rfloor }%
\frac{q^{m\binom{k}{2}}x^{k}y^{n-mk}}{([k]_{q^{m}}!)^{2}[n-mk]_{q}!},\:m,n=0,1,2,\cdots. 
\tag{2.16}
\end{equation}%

\item For a fixed $m=0,1,2,\cdots,$ a generating function for the $q$-2DLP is given  by
\begin{equation}
\sum_{n=0}^{\infty }\frac{t^{n}}{[n]_{q}!}\;_{m}L_{n}(x,y|q)=e_{q}(yt)%
\varepsilon _{q}^{m}(-xt^{m})  \tag{2.17}
\end{equation}%
in terms of the Jackson $q$-exponential function (2.9) and
the $(q,m)$-deformed Bessel-Tricomi function (2.13). See Appendix A for the proof. 
\end{enumerate} 
\section{$q$-deformed 2D Laguerre-Gould-Hopper polynomials}

In this section, we introduce a family of $q$-deformed Laguerre-Gould-Hopper
polynomials generalizing both the  above $q$-2DLP and the  \textit{%
\ }$q$-deformed Gould-Hopper generalized Hermite polynomials which were
defined in [7] by 
\begin{equation}
\mathcal{G}_{n}^{m}(x,y|q)=\left[ n\right] _{q}!\sum_{k=0}^{\left\lfloor 
\frac{n}{m}\right\rfloor }\frac{\left( -1\right) ^{k}q^{m\binom{k}{2}%
}x^{n-mk}y^{k}}{[n-mk]_{q}!\ \ \left[ mk\right]_q !!}.  \tag{3.1}
\end{equation}%
 We shall  establish a  connection
formulae for these polynomials.
\medskip\\
\textbf{Definition 3.1.} \textit{For fixed }$m,n,s=0,1,2,\cdots ,$\textit{\
a family of }$q$\textit{-deformed Laguerre-Gould-Hopper polynomials }$q$-\textit{LGHP} \textit{\ is defined by} 
\begin{equation}
_{L}H_{n}^{\left( m,s\right) }(x,y,z|q):=E_{q^{m}}\left(
(D_{x}^{q^{m}})^{-1}(D_{y}^{q})^{m}\right) (\mathcal{G}_{n}^{s}(x,z|q)). 
\tag{3.2}
\end{equation}%
\textit{or equivalently }%
\begin{equation}
_{L}H_{n}^{\left( m,s\right) }(x,y,z|q):=E_{q^{s}}\left(
z(D_{y}^{q})^{s}\right) \left( _{m}L_{n}(x,y|q)\right).   \tag{3.3}
\end{equation}%
\textit{Explicitly,}%
\begin{equation}
_{L}H_{n}^{\left( m,s\right) }(x,y,z|q)=[n]_{q}!\sum_{k=0}^{\left\lfloor 
\frac{n}{s}\right\rfloor }\frac{q^{s\binom{k}{2}}z^{k}{}_{m}L_{n-sk}(x,y|q)}{%
[k]_{q^{s}}![n-sk]_{q}!}.  \tag{3.4}
\end{equation}%
\textbf{Remark 3.1. }For $z=0,$ the $q$-LGHP
reduce to the above $q$-2DLP. That is, 
\begin{equation}
_{L}H_{n}^{\left( m,s\right) }(x,y,0|q)=_{m}L_{n}(x,y|q).  \tag{3.5}
\end{equation}%
\medskip
\textbf{Proposition 3.1.} \textit{For each fixed }$m,s=0,1,2,\cdots ,$%
\textit{\ the generating function for the }$q$-\textit{LGHP}\textit{ in (3.2)-(3.5) is
given by} 
\begin{equation}
\sum_{n=0}^{\infty }\frac{t^{n}}{[n]_{q}!}\;_{L}H_{n}^{\left( m,s\right)
}(x,y,z|q)=e_{q}(yt)E_{q^{s}}\left( zt^{s}\right) \varepsilon
_{q}^{m}(-xt^{m}),  \tag{3.6}
\end{equation}%
\textit{in terms of the exponential functions} (2.9), (2.10) \textit{and the }$(q,m)$\textit{-deformed
Bessel-Tricomi function} (2.13).
\medskip\\
\textbf{Proof. }We start by inserting $\left( 3.4\right) $ in  the l.h.s  of $\left( 3.6\right).$ This gives%
\begin{equation}
\sum_{n=0}^{\infty }\frac{t^{n}}{[n]_{q}!}\;_{L}H_{n}^{\left( m,s\right)
}(x,y,z|q)=\sum_{n=0}^{\infty }t^{n}\left( \sum_{k=0}^{\left\lfloor \frac{n}{%
s}\right\rfloor }\frac{q^{s\binom{k}{2}}z^{k}\,_{m}L_{n-sk}(x,y|q)}{%
[k]_{q^{s}}![n-sk]_{q}!}\right) .  \tag{3.7}
\end{equation}%
The r.h.s of $\left( 3.7\right) $ also reads successively 
\begin{equation}
\sum_{n=0}^{\infty }\sum_{k=0}^{\left\lfloor \frac{n}{s}\right\rfloor }\frac{%
q^{s\binom{k}{2}}\left( zt^{s}\right) ^{k} t^{n-sk}\,_{m}L_{n-sk}(x,y|q)}{%
[k]_{q^{s}}![n-sk]_{q}!}  \tag{3.8}
\end{equation}%
\begin{equation}
=\sum_{n=0}^{\infty }\sum_{k=0}^{\left\lfloor \frac{n}{s}\right\rfloor }%
\frac{q^{s\binom{k}{2}}\left( zt^{s}\right) ^{k}}{[k]_{q^{s}}!}\frac{_{m}L_{n-sk}(x,y|q)t^{n-sk}}{[n-sk]_{q}!}.  \tag{3.9}
\end{equation}%
Now, by applying the following series manipulation 
\begin{equation}
\sum_{n=0}^{\infty }\sum_{k=0}^{\infty }A\left( k,n\right)
=\sum_{n=0}^{\infty }\sum_{k=0}^{\left\lfloor \frac{n}{m}\right\rfloor
}A\left( k,n-mk\right),  \tag{3.10}
\end{equation}%
where $m$ is any positive integer ([SM], p.101), the  double sum in (3.9)
reads%
\begin{equation}
\sum_{k=0}^{\infty }\frac{q^{s\binom{k}{2}}\left( zt^{s}\right) ^{k}}{%
[k]_{q^{s}}!}\sum_{n=0}^{\infty }\frac{t^{n}}{[n]_{q}!}\,_{m}L_{n}(x,y|q). 
\tag{3.11}
\end{equation}%
Recalling $\left( 2.11\right) $ and $\left( 2.17\right), $ we arrive at the
announced result  $\left( 3.6\right) $. $\square $
\bigskip\\
\textbf{Theorem \ 3.1.} \textit{The}\textbf{\ }$q$-\textit{LGHP}\textit{ in (3.2)
satisfy the connection formula}%
\begin{equation}
_{L}H_{k+l}^{\left( m,s\right) }(x,\xi ,z|q)=\sum\limits_{n,r=0}^{k,l}\left[ 
\begin{array}{c}
k \\ 
n%
\end{array}%
\right] _{q}\left[ 
\begin{array}{c}
l \\ 
r%
\end{array}%
\right] _{q}q^{r\left( r-l\right) }\left( \xi\ominus _{q}y\right)
^{n+r}{}_{m}H_{k+l-n-r}^{\left( m,s\right) }(x,y,z|q).  \tag{3.12}
\end{equation}%
\textbf{Proof.} By replacing $t$ by $u\oplus _{q}t$ in $\left( 3.6\right) $
we get 
\begin{equation}
e_{q}(y\left( u\oplus _{q}t\right) )E_{q^{s}}\left( z\left( u\oplus
_{q}t\right) ^{s}\right) \mathcal{E}_{q}^{m}(-x\left( u\oplus _{q}t\right)
^{m})=\sum_{n=0}^{\infty }\frac{\left( u\oplus _{q}t\right) ^{n}}{[n]_{q}!}%
\;_{L}H_{n}^{\left( m,s\right) }(x,y,z|q).  \tag{3.13}
\end{equation}%
We now apply to the r.h.s of $\left( 3.13\right) $ the identity

\begin{equation}
\sum_{j=0}^{+\infty }F(j)\frac{(x\oplus _{q}y)^{j}}{[j]_{q}!}%
=\sum_{j,s=0}^{+\infty }F(j+s)q^{\binom{s}{2}}\frac{x^{j}y^{s}}{%
[j]_{q}![s]_{q}!}.  \tag{3.14}
\end{equation}%
satisfied by the JHC $q$-addition. So that $\left( 3.13\right) $ becomes 
\begin{equation}
\mathcal{E}_{q}^{m}(-x\left( u\oplus _{q}t\right) ^{m})=\frac{1}{%
e_{q}(y\left( u\oplus _{q}t\right) )E_{q^{s}}\left( z\left( u\oplus
_{q}t\right) ^{s}\right) }\sum_{k,l=0}^{\infty }\frac{u^{k}t^{l}}{%
[k]_{q}![l]_{q}!}\;q^{\binom{l}{2}}\,_{L}H_{k+l}^{\left( m,s\right)
}(x,y,z|q).  \tag{3.15}
\end{equation}%
Note that the r.h.s of (3.15) is independent of  variables $y$ and $z$ so that we
can write for any two variables $\xi ,\zeta $ the following equality 
\begin{equation}
\sum_{k,l=0}^{\infty }\frac{u^{k}t^{l}}{[k]_{q}![l]_{q}!}\;q^{\binom{l}{2}%
}\,_{L}H_{k+l}^{\left( m,s\right) }(x,\xi ,\zeta |q)=\lambda
_{u,t}^{q,s}\left( \xi ,\zeta ;y,z\right) \sum_{k,l=0}^{\infty }\frac{%
u^{k}t^{l}}{[k]_{q}![l]_{q}!}\;q^{\binom{l}{2}}\,_{L}H_{k+l}^{\left(
m,s\right) }(x,y,z|q).  \tag{3.16}
\end{equation}%
where 
\begin{equation}
\lambda _{u,t}^{q,s}\left( \xi ,\zeta ;y,z\right) :=\frac{e_{q}(\xi \left(
u\oplus _{q}t\right) )E_{q^{s}}\left( \zeta \left( u\oplus _{q}t\right)
^{s}\right) }{e_{q}(y\left( u\oplus _{q}t\right) )E_{q^{s}}\left( z\left(
u\oplus _{q}t\right) ^{s}\right) }.  \tag{3.17}
\end{equation}%
By using the  rules in (2.12), on can check that the quantity  $\left( 3.17\right) $ also reads
\begin{equation}
\lambda _{u,t}^{q,s}\left( \xi ,\zeta ;y,z\right) =e_{q}(\left( \xi \ominus
_{q}y\right) \left( u\oplus _{q}t\right) )E_{q^{s}}\left( \left( \zeta
\ominus _{q}z\right) \left( u\oplus _{q}t\right) ^{s}\right).  \tag{3.19}
\end{equation}%
On another hand, the r.h.s  of $\left( 3.19\right) $ coincides with the generating function 
\begin{equation}
\sum_{n,r=0}^{\infty }\frac{u^{n}t^{r}}{[n]_{q}![r]_{q}!}\;q^{\binom{r}{2}}%
\mathcal{G}_{n+r}^{\left( m,s\right) }(\xi \ominus _{q}y,\zeta \ominus
_{q}z|q)  \tag{3.20}
\end{equation}%
involving the $q$-deformed Gould-Hopper generalized Hermite polynomials.
Summarizing the above calculations in $\left( 3.16\right) -\left(
3.20\right) $, we arrive at the sum
\begin{equation*}
\sum_{n,r=0}^{\infty }\frac{u^{n}t^{r}}{[n]_{q}![r]_{q}!}\;q^{\binom{r}{2}}%
\mathcal{G}_{n+r}^{\left( s\right) }(\xi \ominus _{q}y,\zeta \ominus
_{q}z|q)\sum_{k,l=0}^{\infty }\frac{u^{k}t^{l}}{[k]_{q}![l]_{q}!}\;q^{\binom{%
l}{2}}\,_{L}H_{k+l}^{\left( m,s\right) }(x,y,z|q)
\end{equation*}%
\begin{equation}
=\sum_{k,l=0}^{\infty }\frac{u^{k}t^{l}}{[k]_{q}![l]_{q}!}\;q^{\binom{l}{2}%
}\,_{L}H_{k+l}^{\left( m,s\right) }(x,\xi ,\zeta |q).  \tag{3.21}
\end{equation}%
Next, applying the series manipulation ([12], p.100):
\begin{equation}
\sum_{p=0}^{\infty }\sum_{s=0}^{\infty }A\left( p,s\right)
=\sum_{p=0}^{\infty }\sum_{s=0}^{p}A\left( s,p-s\right).  \tag{3.22}
\end{equation}%
 to the l.h.s of $\left( 3.21\right) $, we obtain that%
\begin{equation*}
\sum_{k,l=0}^{\infty }\sum_{n,r=0}^{k,l}\frac{u^{k}t^{l}q^{\binom{r}{2}+%
\binom{l-r}{2}}}{[k]_{q}![l]_{q}!}\begin{bmatrix} k\\n \end{bmatrix}_{q}\begin{bmatrix} l\\r  \end{bmatrix}_{q}\;%
\mathcal{G}_{n+r}^{\left( s\right) }(\xi \ominus _{q}y,\zeta \ominus
_{q}z|q)_{L}H_{k+l-n-r}^{\left( m,s\right) }(x,y,z|q)
\end{equation*}%
\begin{equation}
=\sum_{k,l=0}^{\infty }\frac{u^{k}t^{l}}{[k]_{q}![l]_{q}!}\;q^{\binom{l}{2}%
}\,_{L}H_{k+l}^{\left( m,s\right) }(x,\xi ,\zeta |q).  \tag{3.23}
\end{equation}%
By equating terms with $\displaystyle u^{k}t^{l}/[k]_{q}![l]_{q}!$ and using the
simple combinatorial  fact 
\begin{equation}
\binom{r}{2}+\binom{l-r}{2}-\binom{l}{2}=r\left( r-l\right) ,  \tag{3.24}
\end{equation}%
we arrive at the following result%
\begin{equation}
_{L}H_{k+l}^{\left( m,s\right) }(x,\xi ,\zeta
|q)=\sum_{n,r=0}^{k,l}q^{r\left( r-l\right) }\begin{bmatrix} k\\n \end{bmatrix}_{q}\begin{bmatrix} l\\r  \end{bmatrix}_{q}\;\mathcal{G}_{n+r}^{\left( s\right) }(\xi \ominus _{q}y,\zeta \ominus
_{q}z|q)_{L}H_{k+l-n-r}^{\left( m,s\right) }(x,y,z|q).  \tag{3.25}
\end{equation}%
Putting $\zeta =z$ in the last equation we establish the  result in(3.12). $%
\square $
\bigskip\\
\textbf{Theorem 3.2. }\textit{The following summation formula for the
product of }$q$-\textit{LGHP}
\begin{equation*}
_{L}H_{n}^{\left( m,s\right) }(x,\xi ,\zeta |q)_{L}H_{r}^{\left( m,s\right)
}(X,\Omega ,U|q)=\sum_{k,p=0}^{n,r}\begin{bmatrix} k\\n \end{bmatrix}_{q}\begin{bmatrix} l\\r  \end{bmatrix}_{q}\;%
\mathcal{G}_{k}^{\left( s\right) }(\xi \ominus _{q}y,\zeta \ominus _{q}z|q)
\end{equation*}%
\begin{equation}
\times \mathcal{G}_{p}^{\left( s\right) }(\Omega \ominus _{q}Y,U\ominus
_{q}Z|q)_{L}H_{n-k}^{\left( m,s\right) }(x,y,z|q)_{L}H_{r-p}^{\left(
m,s\right) }(X,Y,Z|q).  \tag{3.26}
\end{equation}%
\textit{holds true.}
\medskip\\
\textbf{Proof.} From the generating function $\left( 3.6\right) ,$ we have%
\begin{equation}
e_{q}(yt)E_{q^{s}}\left( zt^{s}\right) \mathcal{E}%
_{q}^{m}(-xt^{m})e_{q}(YT)E_{q^{s}}\left( ZT^{s}\right) \mathcal{E}%
_{q}^{m}(-XT^{m})  \tag{3.27}
\end{equation}%
\begin{equation*}
=\sum_{n,r=0}^{\infty }\,_{L}H_{n}^{\left( m,s\right)
}(x,y,z|q)_{L}H_{r}^{\left( m,s\right) }(X,Y,Z|q)\frac{t^{n}T^{r}}{%
[n]_{q}![r]_{q}!}.
\end{equation*}%
Replacing in $\left( 3.27\right) $ $y$ by $\xi ,$ $z$ by $\zeta ,$ $Y$ by $%
\Omega $ and $Z$ by $U,$ we get

\begin{equation}
e_{q}(\xi t)E_{q^{s}}\left( \zeta t^{s}\right) \mathcal{E}%
_{q}^{m}(-xt^{m})e_{q}(\Omega T)E_{q^{s}}\left( UT^{s}\right) \mathcal{E}%
_{q}^{m}(-XT^{m})  \tag{3.28}
\end{equation}%
\begin{equation*}
=\sum_{n,r=0}^{\infty }\,_{L}H_{n}^{\left( m,s\right) }(x,\xi ,\zeta
|q)_{L}H_{r}^{\left( m,s\right) }(X,\Omega ,U|q)\frac{t^{n}T^{r}}{%
[n]_{q}![r]_{q}!}.
\end{equation*}%
By replacing 
\begin{eqnarray*}
&&\qquad\qquad\qquad\qquad\qquad \mathcal{E}_{q}^{m} (-xt^m) 
\mathcal{E}_{q}^{m} (-XT^m) =\cr
&& \frac{1}{e_{q}(yt)E_{q^s}(zt^s)\,e_{q}(YT)E_{q^s}(ZT^s)}\sum_{n=0}^{\infty}\sum_{r=0}^{\infty}   {}_LH_n^{(m,s)}(x,y,z|q)\,  {}_LH_r^{(m,s)}(X,Y,Z|q)\frac{ t^n\,T^r}{ [n]_{q}!\, [r]_{q}!}
\end{eqnarray*}
in  the l.h.s. of $\left( 3.28\right)$ and using $\left(2.12\right)$
 one gets, after expanding the exponentials in series, the following
\begin{eqnarray*}
 &&
\sum_{n=0}^{\infty}\sum_{r=0}^{\infty}   {}_LH_n^{(m,s)}(x,\xi ,\zeta|q)\,  {}_LH_r^{(m,s)}(X, \Omega, U|q)\frac{ t^n\,T^r}{ [n]_{q}!\, [r]_{q}!}=\cr
&& \;
e_q\Big( (\xi  \ominus_q  y)t\Big)E_{q^s} \Big((\zeta \ominus_q  z)t^s \Big) e_q\Big( (\Omega  \ominus_q  Y)T\Big)E_{q^s} \Big((U \ominus_q  Z)T^s \Big)\cr
&&\times\;\sum_{n=0}^{\infty}\sum_{r=0}^{\infty}   {}_LH_n^{(m,s)}(x,y,z|q)\,  {}_LH_r^{(m,s)}(X,Y,Z|q)\frac{ t^n\,T^r}{ [n]_{q}!\, [r]_{q}!}\cr
&&=\sum_{n,k=0}^{\infty}\,\mathcal{G}_{k}^s(\xi\ominus_q  y, \zeta\ominus z|q) \,{}_LH_{n}^{(m,s)}(x,y,z|q) \frac{ t^{n+k}}{ [n]_{q}!\, [k]_{q}!}  \cr
&&\times\;\sum_{r,p=0}^{\infty} \,\mathcal{G}_{p}^s(\Omega\ominus_q  Y, U\ominus_q  Z|q) \,{}_LH_{r}^{(m,s)}(X,Y,Z|q)\frac{ T^{r+p}}{ [r]_{q}!\, [p]_{q}!}. \qquad\qquad (3.29)
\end{eqnarray*}
Finally, by replacing $n$ by $n-k$ and $r$ by $r-p$ in the r.h.s. of $\left( 3.29\right) $, the proof is completed.   $\square$

\section{Particular cases}
\begin{enumerate}

\item Putting $z=0$ \ in  Theorem 3.1 leads to a connection formula
for the $q$-2DLP as:
\begin{equation}
_{m}L_{k+l}(x,y|q)=\sum\limits_{n,r=0}^{k,l}\left[ 
\begin{array}{c}
k \\ 
n%
\end{array}%
\right] _{q}\left[ 
\begin{array}{c}
l \\ 
r%
\end{array}%
\right] _{q}q^{r\left( r-l\right) }\left( \xi \ominus _{q}y\right)
^{n+r}\,_{m}L_{k+l-n-r}(x,y|q).  \tag{4.1}
\end{equation}%
In particular, 
\begin{equation}
_{m}L_{k}(x,y|q)=\sum\limits_{n=0}^{k}\left[ 
\begin{array}{c}
k \\ 
n%
\end{array}%
\right] _{q}\left( \xi \ominus _{q}y\right) ^{n}\,_{m}L_{k-n}(x,y|q) 
\tag{4.2}
\end{equation}%
and%
\begin{equation}
_{m}L_{l}(x,y|q)=\sum\limits_{r=0}^{l}\left[ 
\begin{array}{c}
l \\ 
r%
\end{array}%
\right] _{q}q^{r\left( r-l\right) }\left( \xi \ominus _{q}y\right)
^{r}\,_{m}L_{l-r}(x,y|q) \tag{4.3}
\end{equation}%
are obtained by setting $l=0$ and $k=0$ in (4.1) respectively.
\item The following    formulae  for the $q$-LGHP  :

\be
{}_LH_{n}^{(m,s)}(x,\xi,\zeta|q)
= \sum_{k=0}^{n}\begin{bmatrix} n\\k \end{bmatrix}_{q} \mathcal{G}_{k}^s(\xi\om y, \zeta\om z|q)\, {}_LH_{n-k}^{(m,s)}(x,y,z|q) \tag{4.4}
\ee
or 
\be
{}_LH_{n}^{(m,s)}(x,\xi,\zeta|q)
= \sum_{k=0}^{n}\begin{bmatrix} n\\k \end{bmatrix}_{q} q^{k(k-n)} \mathcal{G}_{k}^s(\xi\om y, \zeta\om z|q)\, {}_LH_{n-k}^{(m,s)}(x,y,z|q) \tag{4.5}
\ee
and
\be
{}_LH_{n}^{(m,s)}(x,\xi\op y,\zeta|q)
= \sum_{k=0}^{n}\begin{bmatrix} n\\k \end{bmatrix}_{q}  \xi^{k}\, {}_LH_{n-k}^{(m,s)}(x,y,z|q) \tag{4.6}
\ee
or 
\be
{}_LH_{n}^{(m,s)}(x,\xi\op y,\zeta|q)
= \sum_{k=0}^{n}\begin{bmatrix} n\\k \end{bmatrix}_{q}q^{k(k-n)} \xi^{k}\, {}_LH_{n-k}^{(m,s)}(x,y,z|q) \tag{4.7}
\ee
 are valid.
\item The following formula for the product of two $q$-LGHP :
\begin{equation*}   
 _LH_n^{(m,s)}(x,\xi , \zeta\op z|q)\,  {}_LH_r^{(m,s)}(X, \Omega, U\op Z|q)=  
\end{equation*}
\begin{equation*}\sum_{k,p=0}^{n,r}\begin{bmatrix} n\\k \end{bmatrix}_{q}\begin{bmatrix} r\\p \end{bmatrix}_{q}\,(\xi\om y)^k\,(\Omega\om Y)^p\;{}_LH_{n-k}^{(m,s)}(x,y,z|q) \,{}_LH_{r-p}^{(m,s)}(X,Y,Z|q) \tag{4.8}
\end{equation*}
 holds true. The proof is immediate by replacing 
$\zeta$ by $\zeta\op z$ and $U$ by $U\op Z$ in Eq.(3.26). 
\item Taking $l=0,$ in Eq.(3.25)
and replacing $\xi$ by $\xi\op y, \;\zeta$ by $\zeta\op z$ in the resultant equation, we get
\begin{equation*}
 {}_LH_{k}^{(m,s)}(x,\xi\op y,\zeta\op z|q)
= \sum_{n=0}^{k}\begin{bmatrix} k\\n \end{bmatrix}_{q} \mathcal{G}_{n}^s(\xi, \zeta|q)\, {}_LH_{k-n}^{(m,s)}(x,y,z|q) \tag{4.9}
\end{equation*}
which by taking $\zeta=0,$ yields
\begin{equation*}
 {}_LH_{k}^{(m,s)}(x,\xi\op y, z|q)
= \sum_{n=0}^{k}\begin{bmatrix} k\\n \end{bmatrix}_{q}  \xi^n\, {}_LH_{k-n}^{(m,s)}(x,y,z|q). \tag{4.10}
\end{equation*}
Similarly, taking $k=0,$ in Eq.(3.25)
and replacing $\xi$ by $\xi\op y, \;\zeta$ by $\zeta\op z$ in the resultant equation, we get
\begin{equation*}
 {}_LH_{l}^{(m,s)}(x,\xi\op y, \zeta\op z|q)
= \sum_{r=0}^{l}\begin{bmatrix} l\\r \end{bmatrix}_{q}q^{r(r-l)} \mathcal{G}_{r}^s(\xi, \zeta|q)\, {}_LH_{l-r}^{(m,s)}(x,y,z|q),\tag{4.11}
\end{equation*}
which by taking $\zeta=0,$ yields
\begin{equation*}
 {}_LH_{l}^{(m,s)}(x,\xi\op y,z|q)
= \sum_{r=0}^{l}\begin{bmatrix} l\\r \end{bmatrix}_{q}q^{r(r-l)}\, \xi^r\, {}_LH_{l-r}^{(m,s)}(x,y,z|q).\tag{4.12}
\end{equation*}
\item  Taking $z=0=Z$ in 
  Eq.(3.26) and Eq.(4.8), respectively, we get  the following  formulae  for the $q$-2DLP :
\begin{equation*}
 {}_mL_{k+l}(x,\xi|q)
= \sum_{n,r=0}^{k,l}\begin{bmatrix} k\\n \end{bmatrix}_{q}\begin{bmatrix} l\\r \end{bmatrix}_{q}q^{r(r-l)} (\xi\om y)^{n+r}\,  {}_mL_{k+l-n-r}(x,y|q)\tag{4.13}
\end{equation*}
and
\begin{equation*}
 \qquad\qquad\qquad  {}_mL_n(x,\xi|q)\,  {}_mL_r(X, \Omega|q)=  
 \end{equation*}
\begin{equation*}
\sum_{k,p=0}^{n,r}\begin{bmatrix} n\\k \end{bmatrix}_{q}\begin{bmatrix} r\\p \end{bmatrix}_{q}(\xi\om y)^k(\Omega\om Y)^p\;{}_mL_{n-k}(x,y|q) \,{}_mL_{r-p}(X,Y|q). \tag{4.14}
\end{equation*}
\item Taking $x=0$ and $z\rightarrow\dfrac{-z}{[s]_q}$   in  Eq.(3.25),  
we get the following  formulae 
for the  $q$-deformed \smallskip\\Gould-Hopper generalized Hermite polynomials :
\begin{equation*}
\mathcal{G}_{k+l}^s(\xi,\zeta|q)
= \sum_{n,r=0}^{k,l}\begin{bmatrix} k\\n \end{bmatrix}_{q}\begin{bmatrix} l\\r \end{bmatrix}_{q}q^{r(r-l)} \mathcal{G}_{n+r}^s(\xi\om y, \zeta\om \frac{z}{[s]_q}|q)\,\mathcal{G}_{k+l-n-r}^s(y,z|q). \tag{4.15}
\end{equation*}

\item  Taking  $s=2$ and  $x=0$ in (3.25), we get the following connection formulae for the  $q$-deformed 
   Hermite polynomials $H_n(x,y|q)$ [7]:
\begin{equation*}
 H_{k+l}(\xi,z|q)= \sum_{n,r=0}^{k,l}\begin{bmatrix} k\\n \end{bmatrix}_{q}\begin{bmatrix} l\\r \end{bmatrix}_{q}q^{r(r-l)} (\xi\om y)^{n+r}H_{k+l-n-r}(y,z|q), \tag{4.16}
 \end{equation*}
 \begin{equation*}
H_{n}(\xi,z|q)= \sum_{k=0}^{n}\begin{bmatrix} n\\k \end{bmatrix}_{q} (\xi\om y)^{k}H_{n-k}(y,z|q), \tag{4.17}
\end{equation*}
\begin{equation*}
H_{n}(\xi,z|q)= \sum_{k=0}^{n}\begin{bmatrix} n\\k \end{bmatrix}_{q}q^{k(k-n)} (\xi\om y)^{k}H_{n-k}(y,z|q), \tag{4.18}
\end{equation*}
\begin{equation*}
H_{n}(\xi\op y,z|q)= \sum_{k=0}^{n}\begin{bmatrix} n\\k \end{bmatrix}_{q} \xi^{n-k}H_{k}(y,z|q), \tag{4.19}
\end{equation*}
\begin{equation*}
H_{n}(\xi\op y,z|q)= \sum_{k=0}^{n}\begin{bmatrix} n\\k \end{bmatrix}_{q} q^{k(k-n)}\xi^{n-k}H_{k}(y,z|q). \tag{4.20}
\end{equation*}
Again, taking $s=2,\,x=0=X,\;z=\zeta$ and $U=Z,$ in  equation (3.26), we get that
$$
 \qquad\qquad\qquad H_n(\xi ,z|q)\,  H_r( \Omega, Z|q)=  
 $$
 \begin{equation*}
\sum_{k,p=0}^{n,r}\begin{bmatrix} n\\k \end{bmatrix}_{q}\begin{bmatrix} r\\p \end{bmatrix}_{q}\,(\xi\om y)^k(\Omega\om Y)^pH_{n-k}(y,z|q) \,H_{r-p}(Y,Z|q). \tag{4.21}
\end{equation*}
\item  When  $q\to 1$  in (3.25) and (3.26), we obtain the following summation formulae  for the   Laguerre-Gould-Hopper polynomials ${}_LH_{n}^{(m,s)}(x,y,z)$ [7]:
\begin{equation*}
\label{12}
 {}_LH_{k+l}^{(m,s)}(x,\xi,\zeta)
= \sum_{n,r=0}^{k,l}{k\atopwithdelims () n}{l\atopwithdelims() r}g_{n+r}^s(\xi- y, \zeta- z)\, {}_LH_{k+l-n-r}^{(m,s)}(x,y,z), \tag{4.22}
\end{equation*}
and
\begin{equation*}
  {}_LH_n^{(m,s)}(x,\xi ,\zeta)\,  {}_LH_r^{(m,s)}(X, \Omega, U)=\sum_{k,p=0}^{n,r}{n\atopwithdelims() k}{r\atopwithdelims() p}\,g_{k}^s(\xi- y, \zeta- z)  
\end{equation*}
\begin{equation*}
\times\; g_{p}^s(\Omega- Y, U- Z)\;{}_LH_{n-k}^{(m,s)}(x,y,z) \,{}_LH_{r-p}^{(m,s)}(X,Y,Z). \tag{4.23}
\end{equation*}
\medskip\\
Next, taking $\zeta=z$ and $\zeta=z,\; U=Z$ in (4.22) and (4.23) respectively, one obtains
\begin{equation*}
{}_LH_{k+l}^{(m,s)}(x,\xi,z)
= \sum_{n,r=0}^{k,l}{k\atopwithdelims () n}{l\atopwithdelims() r}  (\xi- y)^{n+r}\, {}_LH_{k+l-n-r}^{(m,s)}(x,y, z) \tag{4.24}
\end{equation*}
and
\begin{equation*}
   \qquad \qquad\qquad \qquad {}_LH_n^{(m,s)}(x,\xi ,z)\,  {}_LH_r^{(m,s)}(X, \Omega, Z)=  
   \end{equation*}
  \begin{equation*} 
\sum_{k,p=0}^{n,r}{n\atopwithdelims() k}{r\atopwithdelims() p}\,(\xi- y)^k\,(\Omega- Y)^p\;{}_LH_{n-k}^{(m,s)}(x,y,z) \,{}_LH_{r-p}^{(m,s)}(X,Y,Z). \tag{4.25}
\end{equation*}

\item Taking  $s=2,\,x=0$ and $q\to 1,$  the  $q-$LGHP are reduced to the higher-order Hermite polynomials. That is, 
\begin{equation*}
{}_LH_{n}^{(m,2)}(0, y, z|q)\equiv g_{n}^m(y,z) \tag{4.26}
\end{equation*}
and  summation formulae (3.25) and (3.26) reduce to the ones  in  [7]. That is,
\begin{equation*}
g_{k+l}^m(\xi,y)= \sum_{n,r=0}^{k,l}{k\atopwithdelims() n}{l\atopwithdelims() r} (\xi- x)^{n+r}g_{k+l-n-r}^m(x,y), \tag{4.27}
\end{equation*}
\begin{equation*}
g_{n}^m(\xi,y)= \sum_{k=0}^{n}{n\atopwithdelims() k}  (\xi- x)^{k}g_{n-k}^m(x,y), \tag{4.28}
\end{equation*}
\begin{equation*}
g_{n}^m(\xi+ x,y)= \sum_{k=0}^{n}{n\atopwithdelims() k} \xi^{n-k}g_{k}^m(x,y) \tag{4.29}
\end{equation*}
and
\begin{equation*}
  g_n^m(\xi ,\zeta)\,  g_r^m(\Omega, U)=
 \sum_{k,p=0}^{n,r}{n\atopwithdelims() k}{r\atopwithdelims() p}\,g_{k}^m(\xi- y, \zeta- z)
 \end{equation*}
 \begin{equation*}
\times\;g_{p}^m(\Omega- Y, U- Z)\;g_{n-k}^m(y,z) \,g_{r-p}^m(Y,Z), \tag{4.30}
\end{equation*}
which, by taking $\zeta=z$ and $U=Z,$ yields
\begin{equation*}
   g_n^m(\xi ,z)\,  g_r^m(\Omega, Z)=  
 \sum_{k,p=0}^{n,r}{n\atopwithdelims() k}{r\atopwithdelims() p}\, (\xi- y)^k\,(\Omega- Y)^p\;g_{n-k}^m(y,z) \,g_{r-p}^m(Y,Z).\tag{4.31}
\end{equation*}
\end{enumerate}
\begin{center}
{\large\textbf{Appendix A}} : \textbf{The proof of (2.17)}\medskip\\
\end{center}

To obtain a closed form for the generating function 
\begin{equation*}
G_q(t;x,y;m):=\displaystyle\sum_{n=0}^{\infty}\frac{t^n}{[n]_q!}\;_mL_n(x,y|q) \tag{A.1}
\end{equation*}

 we start by replacing in the r.h.s of (A.1) the $q$-2DLP $_mL_n(x,y|q)$ by its explicit expression (2.16) as follows

\begin{equation}
G_q(t;x,y;m)=\sum_{n=0}^{%
\infty }t^{n}\left( \sum_{k=0}^{\left\lfloor \frac{n}{m}\right\rfloor }%
\frac{q^{m\binom{k}{2}}x^{k}y^{n-mk}}{([k]_{q^{m}}!)^{2}[n-mk]_{q}!}\right) .
\tag{A.2}
\end{equation}%

The right hand side of $\left( A.2\right) $ also reads 
\begin{equation}
\sum_{n=0}^{\infty }\sum_{k=0}^{\left\lfloor \frac{n}{m}\right\rfloor }%
\frac{q^{m\binom{k}{2}}\left( xt^{m}\right) ^{k}\left( yt\right) ^{n-mk}}{%
([k]_{q^{m}}!)^{2}[n-mk]_{q}!},  \tag{A.3}
\end{equation}%
which can be put in form 
\begin{equation}
\sum_{n=0}^{\infty }\sum_{k=0}^{\left\lfloor \frac{n}{m}\right\rfloor }%
\frac{q^{m\binom{k}{2}}\left( xt^{m}\right) ^{k}}{([k]_{q^{m}}!)^{2}}\frac{%
\left( yt\right) ^{n-mk}}{[n-mk]_{q}!}.  \tag{A.4}
\end{equation}%
Using the same manipulation in (3.10). This
enable us to write $\left( A.4\right) $ as%
\begin{equation}
\sum_{n=0}^{\infty }\sum_{k=0}^{\infty }\frac{q^{m\binom{k}{2}}\left(
xt^{m}\right) ^{k}}{([k]_{q^{m}}!)^{2}}\frac{\left( yt\right) ^{n}}{[n]_{q}!},
\tag{A.5}
\end{equation}%
which coincides with%
\begin{equation*}
e_{q}(yt)\varepsilon _{q}^{m}(-xt^{m})=\sum_{n=0}^{\infty }\frac{q^{m\binom{k%
}{2}}\left( xt^{m}\right) ^{k}}{([k]_{q^{m}}!)^{2}}\sum_{k=0}^{\infty }\frac{%
\left( yt\right) ^{n}}{[n]_{q}!} \tag{A.6}
\end{equation*}
This completes the proof. $\square $

\end{document}